\documentclass[12pt,a4paper,reqno]{article}
\newcommand{\guio}[1]{\nobreakdash-\hspace{0pt}#1}

\usepackage{amsmath,amsthm,amssymb,amscd,epsfig,graphicx}

\newtheorem{theo}{Theorem}

\newtheorem{prop}[theo]{Proposition}
\newtheorem{coro}[theo]{Corollary}

\def\C{\mathbb C}
\def\D{\mathbb D}

\def\H{\mathcal H}

\def\diam{\operatorname{diam}}

\def\Lip{\operatorname{Lip}}

\title{Nonremovable sets for H\"older continuous quasiregular mappings in the plane}
\author{Albert Clop
\thanks{The author was supported by projects MTM2004-00519, HF2004-0208, 2005-SGR-00744
}
}

\date{}

\begin{document}

\maketitle

\begin{abstract}
We show that for any dimension $t>2\frac{1+\alpha K}{1+K}$ there exists a compact set $E$ of dimension~$t$ and a function $\alpha$-H\"older continuous on $\C$, which is $K$-quasiregular only on $\C\setminus E$. To do this, we construct an explicit $K$-quasiconformal mapping that gives, by one side, extremal dimension distortion on a Cantor-type set, and by the other, more H\"older continuity than the usual $1/K$. 
\end{abstract}

\section{Introduction}

Let $\alpha\in(0,1)$. A function $f:\C\rightarrow\C$ is said to be locally $\alpha$-H\"older continuous, that is, $f\in\Lip_\alpha(\C)$, if
\begin{equation}\label{lipalfa}
|f(z)-f(w)|\leq C\,|z-w|^\alpha
\end{equation}
whenever $z,w\in\C$ are such that $|z-w|<1$. A set $E\subset\C$ is said to be {\it{removable}} for $\alpha$\guio{H\"older} continuous analytic functions if every function $f\in\Lip_\alpha(\C)$, holomorphic on $\C\setminus E$, is actually an entire function. It turns out that there is a characterization of these sets $E$ in terms of Hausdorff measures. For $\alpha\in(0,1)$, Dol\v zenko \cite{Do} proved that a set $E$ is removable for $\alpha$-H\"older continuous analytic functions if and only if $\H^{1+\alpha}(E)=0$. When $\alpha=1$, we deal with the class of Lipschitz continuous analytic functions. Although the same characterization holds, a more involved argument, due to Uy \cite{U}, is needed to show that sets of positive area are not removable.\\
The same question may be asked in the more general setting of $K$-quasiregular mappings. Given a domain $\Omega\subset\C$ and $K\geq 1$, one says that a mapping $f:\Omega\rightarrow\C$ is $K$-quasiregular in $\Omega$ if $f$ is a $W^{1,2}_{loc}(\Omega)$ solution of the Beltrami equation,
$$\overline\partial f(z)=\mu(z)\,\partial f(z)$$
for almost every $z\in\Omega$, where $\mu$, the Beltrami coefficient, is a measurable function such that $|\mu(z)|\leq\frac{K-1}{K+1}$ at almost every $z\in\Omega$. If $f$ is a homeomorphism, then $f$ is said to be $K$\guio{quasiconformal}. When $\mu=0$, one recovers the classes of analytic functions and conformal mappings on $\Omega$, respectively.\\
\\
It was shown in \cite{C} that if $E$ is a compact set satisfying $\H^d(E)=0$, $d=2\frac{1+\alpha K}{1+K}$, then $E$ is removable for $\alpha$-H\"older continuous $K$-quasiregular mappings. This means that any function $f\in\Lip_\alpha(\C)$, $K$\guio{quasiregular} in $\C\setminus E$, is actually $K$-quasiregular on the whole plane. To look for results in the converse direction, one observes that any compact set $E$ with $\H^{1+\alpha}(E)>0$ is nonremovable for holomorphic functions and, thus, neither for $K$-quasiregular mappings in $\Lip_\alpha$. Hence, we are interested in dimensions between $d$ and $1+\alpha$. In this paper we show that the index $d$ is sharp in the following sense: given $\alpha\in(0,1)$ and $K\geq 1$, for any $t>d$ there exists a compact set $E$ of dimension~$t$, and a function $f\in\Lip_\alpha(\C)$ which is $K$-quasiregular in $\C\setminus E$, and with no $K$-quasiregular extension to $\C$. In other words, we will construct nonremovable sets at any dimension strictly above~$d$. \\
\\
We first have a look at the case $K=1$. Given a compact set $E$ with $\H^{1+\alpha}(E)>0$, by Frostman's Lemma (see for instance \cite[p.112]{M}), there exists a positive Radon measure $\nu$ supported on~$E$, such that $\nu(B(z,r))\leq C\,r^{1+\alpha}$ for any $z\in E$. Thus, the function $h=\frac{1}{\pi z}\ast \nu$ is $\alpha$-H\"older continuous everywhere, holomorphic outside the support of $\nu$ and has no entire extension.\\
Another close situation is found in the limiting case $\alpha=0$, in which $\Lip_\alpha(\C)$ should be replaced by~$BMO(\C)$. Again, when $K=1$, Kaufman \cite{K} characterized the $BMO$ removable compact sets as those with zero length. When $K>1$, it is known (\cite{ACMOU},\cite{AIM}) that sets with $\H^\frac{2}{K+1}(E)=0$ are removable for~$BMO$ $K$-quasiregular mappings. In fact, the appearence of this index $\frac{2}{K+1}$ is not strange. In~\cite{A}, Astala showed that for any $K$-quasiconformal mapping $\phi$ and any compact set $E$,
\begin{equation}\label{dimdist}
\frac{1}{K}\left(\frac{1}{\dim(E)}-\frac{1}{2}\right)\leq \frac{1}{\dim(\phi(E))}-\frac{1}{2}\leq K\left(\frac{1}{\dim(E)}-\frac{1}{2}\right).
\end{equation}
Furthermore, both equalities are always attainable. In particular, sets of dimension $\frac{2}{K+1}$ are $K$\guio{quasiconformally} mapped to sets of dimension at most $1$, which is the critical point for the analytic $BMO$ situation. Hence, from equality at (\ref{dimdist}), there exists for any $t>\frac{2}{K+1}$ a compact set $E$ of dimension $t$ and a $K$-quasiconformal mapping $\phi$ that maps $E$ to a compact set $\phi(E)$ with dimension 
$$t'=\frac{2Kt}{2+(K-1)t}>1.$$
In particular, $\H^1(\phi(E))>0$. As above, one has a positive Radon measure $\nu$ supported on $\phi(E)$, with linear growth, whose Cauchy transform $h=\frac{1}{\pi z}\ast \nu$ is holomorphic on $\C\setminus E$ and has a $BMO(\C)$ extension which is not entire. Now, since $BMO(\C)$ is invariant under quasiconformal changes of variables \cite{R}, the composition $g=h\circ\phi$ is a $BMO(\C)$ $K$-quasiregular mapping on $\C\setminus E$ which has no $K$-quasiregular extension to $\C$. In other words, the set $E$ is not removable for $BMO$ $K$\guio{quasiregular} mappings. This argument shows that the index $\frac{2}{K+1}$ is somewhat critical for the $BMO$ $K$-quasiregular problem.\\
\\
Our plan is to repeat the above argument, but replacing $BMO(\C)$ by $\Lip_\alpha(\C)$. That is, given any dimension $t>2\frac{1+\alpha K}{1+K}$, we will construct a compact set $E$ of dimension $t$ and a $\Lip_\alpha(\C)$ function which is $K$-quasiregular on $\C\setminus E$ but not on $\C$. We will start with a compact set $E$ of dimension $t$ and a $K$-quasiconformal mapping $\phi$ such that $\dim(\phi(E))=t'=\frac{2Kt}{2+(K-1)t}$. Then, we will show that there are $\Lip_\beta(\C)$ functions, for some $\beta>0$, analytic outside of $\phi(E)$, which in turn induce (by composition) $K$-quasiregular functions on $\C\setminus E$ with some global H\"older continuity exponent. This construction will find two obstacles. First, the extremal dimension distortion of sets of dimension~$2\frac{1+\alpha K}{1+K}$ through $K$-quasiconformal mappings is not exactly $1+\alpha$, the critical number in the analytic setting (this was so for $\alpha=0$). Second, the composition of $\beta$-H\"older continuous functions with $K$-quasiconformal mappings is only in $\Lip_{\beta/K}(\C)$, so there is some loose of regularity that might be decisive. To avoid these troubles, we will construct in an explicit way the mapping $\phi$. This concrete construction allows us to show that $\phi$ exhibits an exponent of H\"older continuity given by
$$
\frac{t}{t'}=\frac{1}{K}+\frac{K-1}{2K}t
$$
which is larger than the usual $\frac{1}{K}$ obtained from Mori's Theorem. This regularity will be sufficient for our purposes. On the other hand, if $\dim(E)=t$ and $\dim(\phi(E))=t'$ it is natural to expect $\phi$ to be $\Lip_{t/t'}$.

\section{Extremal Distortion}

In all this section, $D(z,r)$ will denote the open disk of center $z$ and radius $r$. By $\diam(D)$ we mean the diameter of the disk $D$, and $\lambda D$ is the disk concentric with $D$ but such that $\diam(\lambda D)=|\lambda|\,\diam(D)$. By $\D$ we will mean the unit disk, while $Jf$ will denote the jacobian determinant of the function $f$.\\


Recall that a Cantor-type set $C$ of $m$ components is the only compact set which is invariant through a fixed finite family of $m$ similitudes 
$$\aligned
\varphi_j:&\D\rightarrow\D\\
&z\mapsto\varphi_j(z)=a_j + b_j\,z\endaligned$$
with $a_j,b_j\in\C$ for all $j=1,\dots, m$, such that $\varphi_i(\overline{\D})$ are disjoint. In other words, $C\subset\D$ is the only solution to the equation
$$C=\bigcup_{j=1}^m\varphi_j(C).$$ 
Constructively,
$$C=\bigcap_{N=1}^\infty\left(\bigcup_{\ell(J)=N}\varphi_J(\D)\right)$$
where $\varphi_J=\varphi_{j_N}\circ\dots\circ\varphi_{j_1}$ for any chain $J=(j_1,\dots,j_N)$ of length $\ell(J)=N$ of members of~$\{1,\dots,m\}$. The Hausdorff dimension of $C$ is the only solution $d$ to the equation
$$\sum_{j=1}^m|\varphi'_j|^d=1$$
and under the additional assumption $|\varphi_i'|=r$ for all $i$, one easily gets
$$\dim(C)=d=\frac{\log m}{\log\frac{1}{r}}.$$
Under some additional assumptions on $\varphi_i$ (such as having image sets $\varphi_j(\D)$ uniformly distributed in $\D$), we say that $C$ is a regular Cantor set. For them, $0<\H^d(E)<\infty$.\\
\\
One of the main results in \cite{A} is the sharpness of the dimension distortion equation (\ref{dimdist}). To obtain the equality there, the author moved holomorphically a fixed Cantor-type set $F$. This movement defined actually a {\it{holomorphic motion}} on $F$. An interesting extension result, known as the $\lambda$\guio{lemma}, allowed to extend this motion quasiconformally from $F$ to the whole plane. This procedure avoided most of the technicalities and gave the desired result in a surprisingly direct way. However, since we look both for extremal dimension distortion and higher H\"older continuity, we need to construct $\phi$ explicitly. Thus, let $t\in (0,2)$ and $K\geq 1$ be fixed numbers, and denote $t'=\frac{2Kt}{2+(K-1)t}$. As in \cite{A}, we first give a $K$-quasiconformal mapping $\phi$ that maps a regular Cantor set $E$ of dimension $t$, to another regular Cantor set $\phi(E)$ for which $\dim(\phi(E))$ is as close as we want to $t'$. Later, again as in \cite{A}, we will glue a sequence of such mappings in the convenient way.

\begin{prop}\label{prop1}
Given $t\in (0,2)$, $K\geq 1$ and $\varepsilon>0$, there exists a compact $E$ and a 
$K$\guio{quasiconformal} mapping $\phi:\C\rightarrow\C$, with the following properties:
\begin{enumerate}\itemsep=1pt
\item $\phi$ is the identity mapping on $\C\setminus\D$.
\item $E$ is a self-similar Cantor set, constructed with $m=m(\varepsilon)$ similarities.
\item $\dim(\phi(E))\geq t'-\varepsilon$.
\item $J\phi\in L^p_{loc}(\C)$ if and only if $p\leq\frac{K}{K-1}$.
\item $|\phi(z)-\phi(w)|\leq C\,m^{\frac{1}{t}-\frac{1}{t'}}\,|z-w|^\frac{t}{t'}$ whenever $|z-w|<1$.
\end{enumerate}
\end{prop}
\begin{proof}
Our construction follows the scheme in \cite{AIKM}. Thus, we will obtain $\phi$ as a limit of a sequence of $K$-quasiconformal mappings
$$\phi=\lim_{N\to\infty}\phi_N$$
where every $\phi_N$ will act at the $N$-th step of the construction of $E$. More precisely, both $E$ and $\phi(E)$ will be regular Cantor sets associated to two fixed families of similitudes $(\varphi_j)_{j=1,\dots,m}$ and $(\psi_j)_{j=1,\dots,m}$.
At the $N$-th step, $\phi_N$ will map each {\it{generating disk}} of $E$, $\varphi_{j_1\dots j_N}(\D)$, to the corresponding generating disk of the image set, $\psi_{j_1\dots j_N}(\D)$. Since $\phi$ is said both to be $K$-quasiconformal and to give extremal distortion of dimension, we think about using a typical radial stretching,
$$f(z)=z|z|^{\frac{1}{K}-1}$$
conveniently modified. It turns out that this radial stetching $f$ is extremal for some basic properties of $K$-quasiconformal mappings, such as H\"older continuity. In order to avoid that our mapping $\phi$ also lives in no better H\"older space than $\Lip_{1/K}(\C)$ (actually for $f$ this is the case), we will replace~$f$ by a linear mapping in a small neighbourhood of its singularity. This change will not affect the integrability index, and at the same time will give some improvement on the H\"older exponent.\\

\vspace*{5pt}
\noindent
Take $m\geq 100$, and consider $m$ disjoint disks inside of $\D$, $D(z_i, r)$, uniformly distributed, all with the same radius $r=r_m$. By taking $m$ big enough, we may always assume that $c_m=mr^2\geq\frac{1}{2}$.
Given any $\sigma\in (0,1)$ to be determined later, we can consider $m$ similitudes
$$\varphi_i(z)=z_i+\sigma r\,z,\,\,\,\,\,\,\,z\in\D$$
and denote, for every $i=1,\dots, m$,  
$$\aligned
D_i&=\frac{1}{\sigma}\,\varphi_i(\D)=D(z_i, r_1)\\
D'_i&=\varphi_i(\D)=D(z_i,\sigma r_1)
\endaligned$$
where here we have written $r_1=r$. We define
$$g_1(z)=\begin{cases}
\sigma^{\frac{1}{K}-1}(z-z_i)+z_i&z\in D'_i\\
\left|\frac{z-z_i}{r_1}\right|^{\frac{1}{K}-1}(z-z_i)+z_i&z\in D_i\setminus D'_i\\
z&\text{otherwise}.
\end{cases}$$
It may be easily seen that $g_1$ defines a $K$-quasiconformal mapping, which is conformal everywhere except on each ring $D_i\setminus D'_i$. 

\begin{figure}[ht]
\begin{center}
\includegraphics{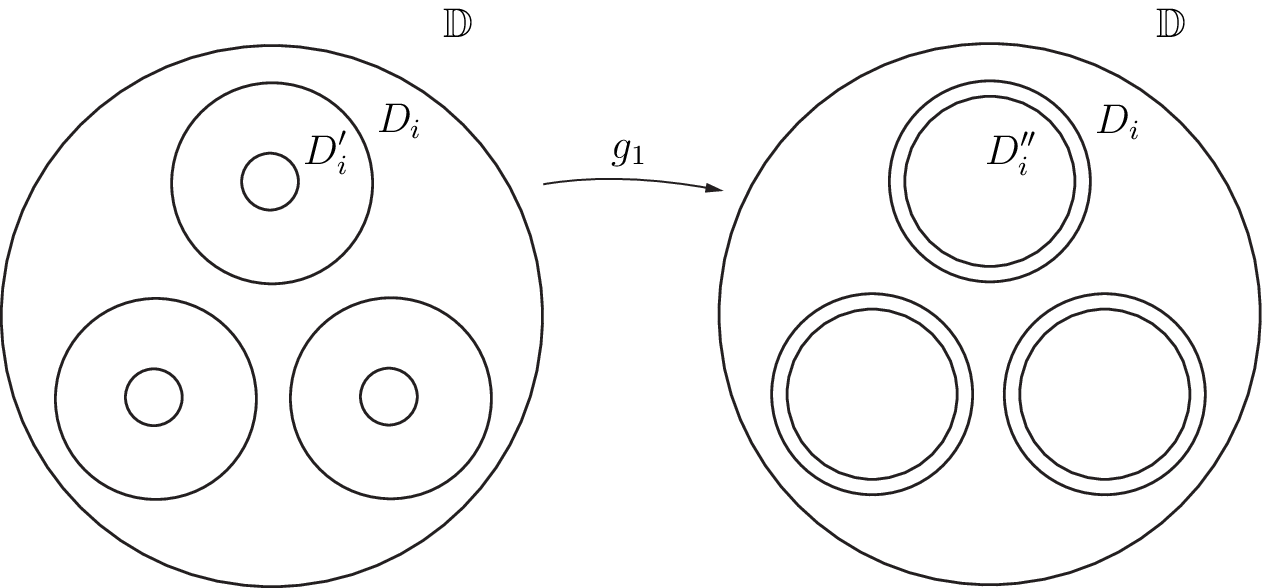}\\
Figure 1
\end{center}
\end{figure}

\noindent
Moreover, if we put
$$\psi_i(z)=z_i+\sigma^\frac{1}{K}r\,z,\,\,\,\,\,\,z\in\D.$$
then $g_1$ maps every $D_i$ to itself, while each $D'_i$ is mapped to~$D''_i=\psi_i(\D)$, as is shown in Figure 1. Now we denote $\phi_1=g_1$. \\
\\
At the second step, we will repeat the above procedure inside of every $D''_i$, and the rest will remain fixed. That is, we will define $g_2$ on the target set of $\phi_1$, and then construct $\phi_2$ as
$$\phi_2=g_2\circ\phi_1$$
To do it more explicitly, we deonte
$$\aligned
D_{ij}&=\frac{1}{\sigma}\,\phi_1\left(\varphi_{ij}(\D)\right)=D(z_{ij}, r_2)\\
D'_{ij}&=\phi_1\left(\varphi_{ij}(\D)\right)=D(z_{ij},\sigma r_2)
\endaligned$$
where a computation shows that $r_2=\sigma^\frac{1}{K}r\,r_1$. Now we define 
$$g_2(z)=\begin{cases}
\sigma^{\frac{1}{K}-1}(z-z_{ij})+z_{ij}&z\in D'_{ij}\\
\left|\frac{z-z_{ij}}{r_2}\right|^{\frac{1}{K}-1}(z-z_{ij})+z_{ij}&z\in D_{ij}\setminus D'_{ij}\\
z&\text{otherwise}.
\end{cases}$$
By construction, $g_2$ is $K$-quasiconformal on $\C$, conformal outside a union of $m^2$ rings, and maps~$D'_{ij}$ to~$D''_{ij}=\psi_{ij}(\D)$, while every point outside of the disks $D_{ij}$ remains fixed through $g_2$, as Figure 2 shows. 

\begin{figure}[ht]
\begin{center}
\includegraphics{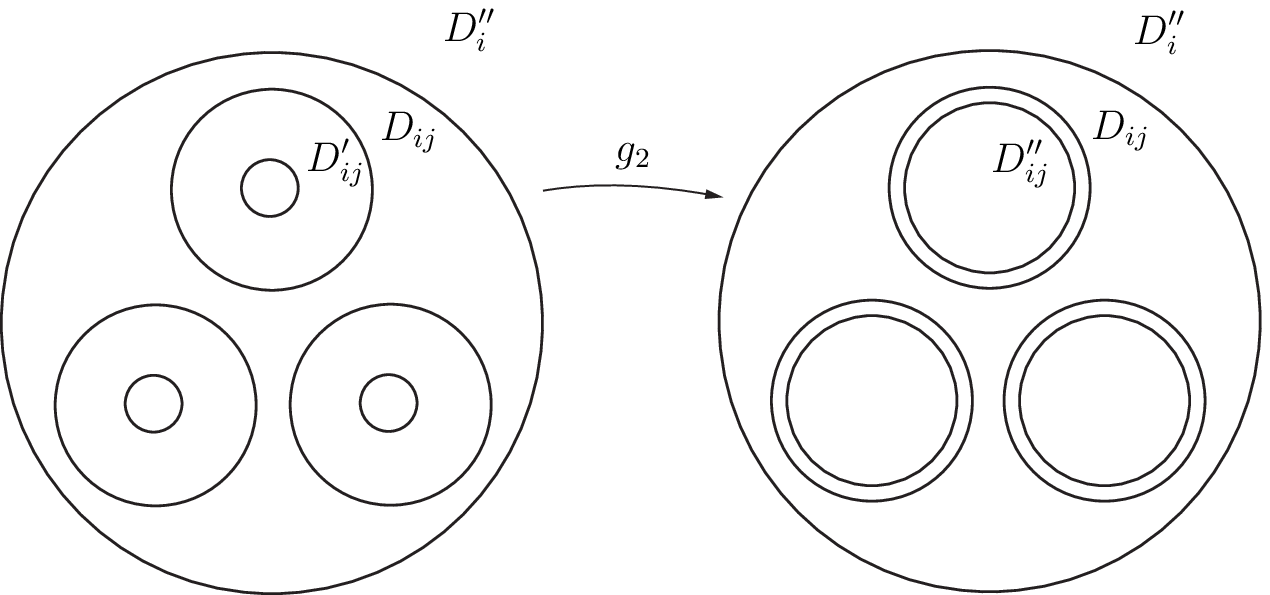}\\
Figure 2
\end{center}
\end{figure}

\noindent
Thus, the composition $\phi_2=g_2\circ\phi_1$ (see Figure 3) is still $K$\guio{quasiconformal}, agrees with the identity outside of $\D$, and
$$\phi_2(\varphi_{ij}(\D))=\psi_{ij}(\D)$$
for any $i,j=1,\dots,m$.\\
\\
After $N-1$ steps, we will define $g_N$ on the target side of $\phi_{N-1}$. For each multiindex $J=(j_1,...,j_N)$ of length $\ell(J)=N$ we will denote
$$\aligned
D_{J}&=\frac{1}{\sigma}\,\phi_{N-1}\left(\varphi_{J}(\D)\right)=D(z_{J}, r_N)\\
D'_{J}&=\phi_{N-1}\left(\varphi_{J}(\D)\right)=D(z_{J},\sigma r_N)
\endaligned$$
where now $r_N=\sigma^\frac{1}{K}r\,r_{N-1}$. Then, the mapping
$$g_N(z)=\begin{cases}
\sigma^{\frac{1}{K}-1}(z-z_{J})+z_{J}&z\in D'_{J}\\
\left|\frac{z-z_{J}}{r_N}\right|^{\frac{1}{K}-1}(z-z_{J})+z_{J}&z\in D_{J}\setminus D'_{J}\\
z&\text{otherwise}
\end{cases}$$
is $K$-quasiconformal on the plane, with $K>1$ only on a union of $m^N$ rings. Moreover, $g_N(D_J)=D_J$ and $g_N(D'_J)=D''_J$, where $D''_J=\psi_J(\D)$.
 
\begin{figure}[ht]
\begin{center}
\includegraphics{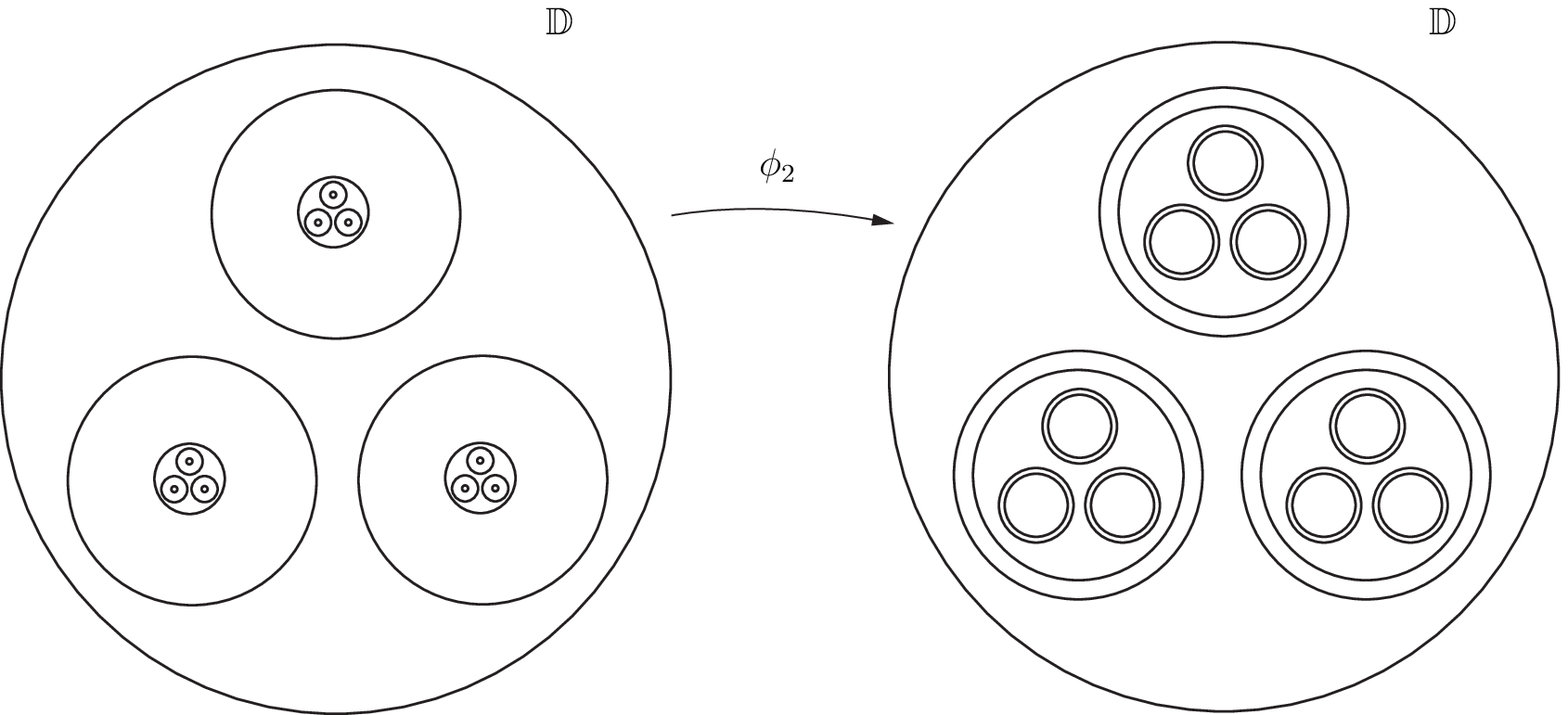}\\
Figure 3
\end{center}
\end{figure}

\noindent
 
As a consequence, the composition $\phi_N=g_N\circ\phi_{N-1}$ is also $K$-quasiconformal and
$$\phi_N(\varphi_J(\D))=\psi_J(\D).$$
With this procedure, it is clear that the sequence $\phi_N$ is uniformly convergent to a homeomorphism~$\phi$. It is also clear that $\phi$ has distortion bounded by $K$ almost everywhere and, in fact, that $\phi$ is a $K$-quasiconformal mapping. By construction, $\phi$ maps the regular Cantor set
\begin{gather*}
E=\bigcap_{N=1}^\infty\left(\bigcup_{\ell(J)=N}\varphi_J(\D)\right)\\
\intertext{to}
\phi(E)=\bigcap_{N=1}^\infty\left(\bigcup_{\ell(J)=N}\psi_{J}(\D)\right)
\end{gather*}
which obviously is also a regular Cantor set. If now we choose $\sigma$ so that
$$m(\sigma r)^t=1$$
we directly obtain, by one hand, $0<\H^t(E)<\infty$, and by the other,
$$
\frac{1}{\dim(\phi(E))}=\frac{1}{t'}+\frac{K-1}{2K}\,\frac{\log\frac{1}{mr^2}}{\log m}.
$$
Since $c_m\geq\frac{1}{2}$ for all $m$, we may always get
$$\dim(\phi(E))\geq t'-\varepsilon$$
just increasing $m$ if needed.\\
\\


Now we have to look at the regularity properties of our mapping $\phi$. To do that, we introduce the following notation: put $G^0=\D$, and denote by $P^N_J$ and by $G^N_J$, respectively, the protecting and generating disks of generation $N$, that is, for any chain $J=(j_1,...,j_N)$,
$$\aligned
P^N_J&=\frac{1}{\sigma}\,\varphi_J(\D)\\G^N_J&=\varphi_J(\D).
\endaligned$$
With this notation, $D_J=\phi_{N-1}(P^N_J)$, $D'_J=\phi_{N-1}(G^N_J)$ and $D''_J=\phi_N(G^N_J)$.\\
\\
Now take any $p$ such that $J\phi\in L^p_{loc}(\C)$. Of course, we can assume $p\geq 1$. Then, one may decompose the $p$-mass of $J\phi$ over $\D$ in the following way
$$
\int_\D J\phi(z)^pdA(z)=\int_{\D\setminus\cup_iP^1_i}J\phi(z)^pdA(z)+\sum_{i=1}^m\int_{P^1_i\setminus G^1_i}J\phi(z)^pdA(z)+\sum_{i=1}^m\int_{G^1_i}J\phi(z)^pdA(z)\\
$$
and since $\phi=\phi_1$ on $\C\setminus\cup_iG^1_i$,  
$$
\int_\D J\phi(z)^pdA(z)=\int_{\D\setminus\cup_iP^1_i}J\phi_1(z)^pdA(z)+m\int_{P^1\setminus G^1}J\phi_1(z)^pdA(z)+m\int_{G^1}J\phi(z)^pdA(z)
$$
where here $P^1$ and $G^1$ denote, respectively, any of the protecting and generating disks of first generation. Now, one may repeat this argument for the last integral, and by a recursive argument we get
$$
\int_\D J\phi(z)^pdA(z)=\sum_{N=0}^\infty\,m^N\,\int_{G^N\setminus\cup_iP^{N+1}_i}J\phi_{N+1}(z)^pdA(z)+\sum_{N=1}^\infty\,m^N\,\int_{P^N\setminus G^N}J\phi_N(z)^pdA(z)
$$
where, as before, $P^N$ and $G^N$ denote any protecting or, respectively, generating disk of $N$-th~generation. Now we compute separatedly the integrals in both sums. By one hand, if $J=(j_1,...,j_N)$
$$\aligned
\int_{P^N_J\setminus G^N_J}J\phi_N(z)^pdA(z)
&=\int_{P^N_J\setminus G^N_J}Jg_N(\phi_{N-1}(z))^p\,J\phi_{N-1}(z)^pdA(z)\\
&=\int_{D_J\setminus D'_J}Jg_N(w)^p\,J\phi_{N-1}(\phi_{N-1}^{-1}(w))^{p-1}dA(w)\\
&=\left(\sigma^{\frac{1}{K}-1}\right)^{2(N-1)(p-1)}\int_{D_J\setminus D'_J}Jg_N(w)^pdA(w)\\
&=r^{2N}\,\sigma^{(N-1)\gamma}\,\,\frac{2\pi}{K^p}\,\left|\frac{1-\sigma^{\gamma}}{\gamma}\right|\\
\endaligned$$
under the additional assumption $p\neq\frac{K}{K-1}$, and where $\gamma=2p\left(\frac{1}{K}-1\right)+2$. If $p=\frac{K}{K-1}$, we get
$$
\int_{P^N_J\setminus G^N_J}J\phi_N(z)^\frac{K}{K-1}dA(z)=r^{2N}\,\frac{2\pi}{K^\frac{K}{K-1}}\,\log\frac{1}{\sigma}.
$$
On the other hand, for any value of $p$,
$$\aligned
\int_{G^N_J\setminus\cup_iP^{N+1}_{(Ji)}}J\phi_{N+1}(z)^pdA(z)
&=\int_{G^N_J\setminus\cup_iP^{N+1}_{(Ji)}}Jg_{N+1}(\phi_N(z))^p\,J\phi_N(z)^pdA(z)\\
&=\int_{D''_J\setminus\cup_iD_{(Ji)}}Jg_{N+1}(w)^p\,J\phi_N(\phi_N^{-1}(w))^{p-1}dA(w)\\
&=\left(\sigma^{\frac{1}{K}-1}\right)^{2N(p-1)}\int_{D''_J\setminus\cup_iD_{(Ji)}}Jg_{N+1}(w)^pdA(w)\\
&=\left(\sigma^{\frac{1}{K}-1}\right)^{2N(p-1)}\int_{D''_J\setminus\cup_iD_{(Ji)}}1\,dA(w)\\
&=\left(\sigma^{\frac{1}{K}-1}\right)^{2N(p-1)}\,|D''_J\setminus\cup_iD_{(Ji)}|\\
&=r^{2N}\,\sigma^{N\gamma}\,\,\pi(1-c_m).
\endaligned$$
Thus, for any $p\neq\frac{K}{K-1}$ we get
$$
\int_\D J\phi(z)^pdA(z)=\left(\pi(1-c_m)+c_m\,\frac{2\pi}{K^p}\,\left|\frac{1-\sigma^\gamma}{\gamma}\right|\right)\,\sum_{N=0}^\infty\left(c_m\,\sigma^\gamma\right)^N.
$$
Since $p$ is such that $J\phi\in L^p_{loc}(\C)$, we necessarily get $\sigma^\gamma<\frac{1}{c_m}$. For $m$ big enough this is equivalent to $\gamma>0$, i.e. $p<\frac{K}{K-1}$. At the critical point, $p=\frac{K}{K-1}$, one gets
$$
\int_\D J\phi(z)^\frac{K}{K-1}dA(z)=\left(\pi(1-c_m)+c_m\,\frac{2\pi}{K^\frac{K}{K-1}}\,\log\frac{1}{\sigma}\right)\,\sum_{N=0}^\infty\left(c_m\right)^N
$$
which will always converge, for any fixed value of $m$. This shows that we can choose $m$ big enough so that $J\phi\in L^p_{loc}(\D)$ if and only if $p\leq\frac{K}{K-1}$.\\
\\
Finally, it just remains to check that $\phi$ is H\"older continuous with exponent $\gamma=t/t'$. By means of Poincar\'e inequality together with the quasiconformality of $\phi$, it is enough \cite[p.64]{G} to show that for any disk $D$
$$\int_DJ\phi(z)\,dA(z)\leq C\,\diam(D)^{2t/t'}.$$
Hence, for some fixed disk $D$, take $N$ such that $(\sigma r)^N\leq\frac{1}{2}\diam(D)<(\sigma r)^{N-1}$. We have
$$
\int_D J\phi(z)\,dA(z)\leq\int_{D\setminus\cup G^N_J}J\phi(z)\,dA(z)+\int_{\cup G^N_J}J\phi(z)\,dA(z)
$$
where the union $\cup G^N_J$ runs over all disks $G^N_J$ such that $G^N_J\cap D\neq\emptyset$. On $D\setminus\cup G^N_J$, we easily see that
$$J\phi=J\phi_N\leq\frac{1}{K}\,\left(\sigma^{\frac{1}{K}-1}\right)^{2N}.$$
Thus,
$$\aligned
\int_{D\setminus\cup G^N_J}J\phi(z)\,dA(z)
&\leq\frac{1}{K}\left(\sigma^{\frac{1}{K}-1}\right)^{2N}\,\pi\left(\frac{1}{2}\diam(D)\right)^2\\
&\leq\frac{\pi}{K}\left(c_m^\frac{K-1}{2K}\right)^{2N}\,m^{2(\frac{1}{t}-\frac{1}{t'})}\,\left(\frac{1}{2}\diam(D)\right)^{2\frac{t}{t'}}.
\endaligned$$
On the other hand, recall that $\phi(G^N_J)=\phi_N(G^N_J)$ are disks of radius $\left(\sigma^\frac{1}{K}r\right)^N$. Hence,
$$\aligned
\int_{\cup_JG^N_J}J\phi(z)\,dA(z)&=\sum_J\int_{G^N_J}J\phi(z)\,dA(z)=\sum_J|\phi(G^N_J)|\\
&=\sum_J|\phi_N(G^N_J)|=\sum_J\pi\left(\sigma^\frac{1}{K}r\right)^{2N}\\
&=\pi\left(c_m^\frac{K-1}{2K}\right)^{2N}\,\left(\frac{1}{2}\diam(D)\right)^{t/t'}\sum_J\left(\frac{(\sigma r)^N}{\frac{1}{2}\diam(D)}\right)^{2t/t'}\,\\
\endaligned$$
and it just remains to bound $\sum_J\left(\frac{(\sigma r)^N}{\frac{1}{2}\diam(D)}\right)^{2t/t'}$. Actually, this is equivalent to find some constant~$C$ such that
$$
\sum_{G^N_J\cap D\neq\emptyset}\diam(G^N_J)^{2t/t'}\leq C\,\diam(D)^{2t/t'}.
$$
But the disks $G^N_J$ come from a self-similar construction, said to give a regular Cantor set of dimension $t$. In particular, they may be chosen uniformly distributed so that the so-called $t$-{\it{dimensional packing condition}} is satisfied, that is,
$$
\sum_{G^N_J\cap D\neq\emptyset}\diam(G^N_J)^{t}\leq C\,\diam(D)^{t}.
$$
It is easy to show that this condition implies the $s$-dimensional one for all $s>t$ (in particular, for~$s=\frac{2t}{t'}$). Hence, the constant $C$ exists and is independent of $m$. Thus, what we finally get is that
$$\int_DJ\phi(z)\,dA(z)\leq C\,m^{\frac{1}{t}-\frac{1}{t'}}\,\left(c_m^\frac{K-1}{2K}\right)^{2N}\,\left(\frac{1}{2}\diam(D)\right)^{2t/t'}$$
and the result follows.
\end{proof}

\newpage

\begin{coro}\label{coro1}
Let $K\geq 1$ and $t\in(0,2)$. There exists a compact set $E$, of dimension $t$, and a $K$-quasiconformal mapping $\phi:\C\rightarrow\C$, such that:
\begin{enumerate}
\item $\H^t(E)$ is $\sigma$-finite.
\item $\dim(\phi(E))=t'$.
\item $|\phi(z)-\phi(w)|\leq C\,|z-w|^\frac{t}{t'}$ whenever $|z-w|<1$.
\end{enumerate}
\end{coro}
\begin{proof}
Given $\varepsilon>0$, $K\geq 1$ and $t\in(0,2)$, let $\phi:\C\rightarrow\C$ and $E$ be as in Proposition \ref{prop1}. Then, for any fixed $r>0$, the mapping
$$\psi_r(z)=r\,\phi(z/r)$$
and the set $E_r=rE$ exhibit the same properties than $\phi$ and $E$, since neither $K$-quasiconformality nor Hausdorff dimension are modified through dilations. However, when computing the new $\Lip_{t/t'}$ constant, if $|z-w|<r$ then
$$
|\psi_r(z)-\psi_r(w)|=r\,|\phi(z/r)-\phi(w/r)|\leq C\, m^{\frac{1}{t}-\frac{1}{t'}}\,r^{1-\frac{t}{t'}}\,|z-w|^\frac{t}{t'}.
$$
Thus, as in \cite{A}, let $D_j=D(z_j,r_j)$ be a countable disjoint family of disks inside of $\D$, and let $\varepsilon_j$ be a sequence of positive numbers, $\varepsilon_j\to 0$ as $j\to\infty$. For each $j$, let $\phi_j$ and $E_j$ be as in Proposition~\ref{prop1}, so that $\dim(\phi_j(E_j))\geq t'-\varepsilon_j$. In particular, each $E_j$ is a regular Cantor set of $m_j$ components. Denote then $\psi_j(z)=r_j\,\phi_j(\frac{z-z_j}{r_j})$ and $F_j=z_j+r_j\,E_j$, and define
$$\psi(z)=\begin{cases}
\psi_j(z)&z\in D_j\\
z&\text{otherwise}.\end{cases}$$
By construction, $\psi$ is a $K$-quasiconformal mapping. It maps the set $F=\cup_jF_j$ to the set $\psi(F)=\cup_j\psi_j(F_j)$. Moreover, $\H^t(F)$ is $\sigma$-finite, while
$$\dim(\phi(F))=\sup_j\,\dim(\psi_j(F_j))=t'.$$
Finally, assume that $z$ lives inside some fixed $D_k$ and that $w\in\D\setminus\cup_jD_j$. Then, consider the line segment $L$ between $z$ and $w$, and denote $\{z_k\}=L\cap\partial D_k$. Then, both $z_k$ and $w$ are fixed points for $\psi$, so that
$$\aligned
|\psi(z)-\psi(w)|&\leq|\psi(z)-\psi(z_k)|+|\psi(z_k)-\psi(w)|\\
&\leq C\,m_k^{\frac{1}{t}-\frac{1}{t'}}\,r_k^{1-\frac{t}{t'}}\,|z-z_k|^\frac{t}{t'}+|z_k-w|.
\endaligned$$
Since we are still free of choosing the disks $D_j$, we may do it so that the radius $r_j$ satisfy
$$m_k^{\frac{1}{t}-\frac{1}{t'}}\,r_k^{1-\frac{t}{t'}}<1$$
or, equivalently, $m_j\,r_j^t<1$. Under this assumption, we finally get
$$|\psi(z)-\psi(w)|\leq (C+1)\,|z-w|^\frac{t}{t'}$$
whenever $|z-w|<1$. This clearly shows that $\psi\in\Lip_{t/t'}(\C)$.
\end{proof}

Although the set in Corollary \ref{coro1} is more critical than the one we constructed in Proposition \ref{prop1}, in the sense that the first gives precisely the extremal dimension distortion, both do the same work when studying non-removable sets for H\"older-continuous quasiregular mappings.

\begin{coro}
Let $K\geq 1$ and $\alpha\in(0,1)$. For any $t>2\frac{1+\alpha K}{1+K}$ there exists a compact set $E$ with~$0<\H^t(E)<\infty$, non removable for $K$-quasiregular mappings in $\Lip_\alpha$.
\end{coro}
\begin{proof}
Let $E$ and $\phi$ be such that $\dim\phi(E)\geq t'-\varepsilon>1$ for some $\varepsilon$ small enough. Hence, by Frostmann's Lemma, we can construct a positive Radon measure $\mu$ supported on $\phi(E)$, with growth~$t'-2\varepsilon$. Its Cauchy transform $g={\cal C}\mu$ defines a holomorphic function on $\C\setminus\phi(E)$, not entire, and with a H\"older continuous extension to the whole plane, with exponent $t'-2\varepsilon -1$. Set
$$f=g\circ\phi.$$
Clearly, $f$ is $K$-quasiregular on $\C\setminus E$ and has no $K$-quasiregular extension to $\C$. Furthermore, $f$ is H\"older continuous with exponent
$$(t'-2\varepsilon -1)\frac{t}{t'}=t-(2\varepsilon+1)\frac{t}{t'}.$$
Thus, we just need $\varepsilon>0$ small enough so that
$$
t-(2\varepsilon+1)\frac{t}{t'}\geq\alpha
$$
but this inequality is equivalent to
$$
\left(t-2\frac{1+\alpha K}{1+K}\right)\geq\varepsilon\,\frac{2}{K+1}\,\left(2+(K-1)t\right)
$$
and the proof is complete.
\end{proof}

Something similar may be said when dealing with finite distortion mappings. Recall that if $\Omega\subset\C$ is an open set, then a {\it{finite distortion mapping}} on $\Omega$ is a function $f:\Omega\rightarrow\C$ in the Sobolev class $W^{1,1}_{loc}(\C)$ with locally integrable jacobian, $Jf\in L^1_{loc}(\C)$, and such that there is a measurable function $K_f:\Omega\rightarrow[1,\infty]$, finite almost everywhere, called the {\it{distortion function}} of $f$, for which
$$|Df(z)|^2\leq K_f(z)\,Jf(z)$$
at almost every $z\in\Omega$. When $K_f\in L^\infty$, $\|K_f\|_\infty=K$, one recovers the class of $K$-quasiregular mappings. However, weaker assumptions on $K_f$ also give interesting results. The most typical situation appears for {\it{exponentially integrable distortion mappings}}, that is, finite distortion mappings~$f$ for which the distortion function $K_f$ satisfies
$$\exp\{K_f\}\in L^p_{loc}(\C)$$
for some $p$ big enough. In \cite{C}, it was shown that compact sets $E$ with $\sigma$-finite $\H^{2\alpha}(E)$ are removable for exponentially integrable distortion mappings in $\Lip_\alpha$. 

\begin{coro}
Let $\alpha\in (0,1)$. For any $t>2\alpha$ there exists a compact set $E$ of dimension $t$ and a function $f\in\Lip_\alpha(\C)$, which defines an exponentially integrable distortion mapping on $\C\setminus E$, and has no finite distortion extension to $\C$.
\end{coro}
\begin{proof}
If $t>2\alpha$, then there exists $K\geq 1$ such that $t>2\frac{1+\alpha K}{1+K}$. Thus, we have a compact set $E$ of dimension $t$, and a $\Lip_\alpha(\C)$ function $f$, $K$-quasiregular on $\C\setminus E$ but not on $\C$. Of course, $f$ is an exponentially integrable distortion mapping on $\C\setminus E$, with distortion function $K_f$ essentially bounded by $K$. If $f$ extended to a finite distortion mapping on $\C$, in particular we would have $Jf\in L^1_{loc}(\C)$. But then, since $K_f\leq K$ at almost every point, this would imply that actually $f$~extends $K$-quasiregularly.
\end{proof}

At this point, it should be said that above the critical index $2\frac{1+\alpha K}{1+K}$ one might find also some removable set. For instance, due to an unpublished result of S. Smirnov, it is known that if $E=\partial\D$ and $\phi$ is a $K$-quasiconformal mapping, then
$$\dim(\phi(E))\leq 1+\left(\frac{K-1}{K+1}\right)^2$$
which is better than the usual dimension distortion equation (\ref{dimdist}). Hence, if we choose  $K\geq 1$ small enough, then there exists $\alpha$ satisfying
$$K\left(\frac{K-1}{K+1}\right)^2<\alpha<\frac{K-1}{2K}.$$
For those values of $\alpha$, the set $E=\partial\D$ is removable for $\alpha$-H\"older continuous $K$-quasiregular mappings and, however, 
$$2\frac{1+\alpha K}{1+K}<\dim(E).$$
This suggests that between $2\frac{1+\alpha K}{1+K}$ and $1+\alpha$ everything may happen. \\
\\
{\bf{Aknowledgements}}. Part of this work was done while the author was visiting the mathematics departments at the universities of Helsinki and Jyv\"askyl\"a (Finland). Thanks are due to both institutions for their hospitality. The author is also grateful to their advisors J. Mateu and J. Orobitg, as well as to D. Faraco and X. Zhong, for many interesting discussions on the subject of the paper.

\vskip 1cm
\begin{itemize}

\item[]{Departament de Matem\`atiques, Facultat de Ci\`encies\\Campus de la Universitat Aut\`onoma de Barcelona\\ 08193-Bellaterra, Barcelona (Spain)\\
 {albertcp@mat.uab.es}\\ Tel. +34 93 581 45 45\\Fax +34 93 581 27 90}

\end{itemize}
\end{document}